\documentclass[a4paper,11pt]{article}
\title{Matrices related to the Pascal triangle}
\author{Roland Bacher}
\begin{document}
\maketitle
%\centerline{Revised version of
%pr\'epublication de l'Institut Fourier $\hbox{n}^0$ 541 (2001)}
%\centerline{\it http://www-fourier.ujf-grenoble.fr/prepublications.html}  

{\it Abstract\footnote{Math. Class: 11B39, 11B65, 11C20 
}: 
The aim of this paper is to study determinants of
matrices related to the Pascal triangle.}
\footnote
{Support from the Swiss National Science Foundation is gratefully acknowledged.
}

\section{The Pascal triangle}

Let $P$ be the  infinite symmetric \lq\lq matrix'' with entries 
$p_{i,j}={i+j\choose i}$ for $0\leq i,j\in {\bf N}$. The matrix $P$ is hence the
famous Pascal triangle yielding the binomial coefficients and can be recursively
constructed by the rules $p_{0,i}=p_{i,0}=1$ for $i\geq 0$ and
$p_{i,j}=p_{i-1,j}+p_{i,j-1}$ for $1\leq i,j$. 

In this paper we are interested in (sequences of determinants of finite) matrices
related to $P$.

The present section deals with some minors (determinants of submatrices)
of the above Pascal triangle $P$, perhaps slightly perturbed.

Sections 2-6 are devoted to the study of matrices satisfying the Pascal
recursion rule $m_{i,j}=m_{i-1,j}+m_{i,j-1}$ for $1\leq i,j<n$ (with
various choices for the first row $m_{0,j}$ and column $m_{i,0}$). Our main
result is the experimental observation (Conjecture 3.3 and Remarks 3.4)
that given such an infinite matrix
whose first row and column satisfy linear recursions (like for instance the
Fibonacci sequence $1,1,2,3,5,8,13,21,\dots$), 
then the determinants of a suitable sequence of 
submatrices seem also to satisfy a linear recursion. 
We give a proof if all linear recursions 
are of length at most $2$ (Theorem 3.1).

Section 7 is seemingly unrelated since it deals with matrices which are
\lq\lq periodic'' along strips parallel to the diagonal. If such a matrix consists only of 
a finite number of such strips, then an appropriate sequence of determinants
satisfies a linear recursion (Theorem 7.1). 

Section 8 is an application of section 7. It deals with matrices which are periodic
on the diagonal and off-diagonal coefficients satisfy a different kind of
Pascal-like relation.

We come now back to the Pascal triangle $P$ with coefficients
$p_{i,j}={i+j\choose i}$. Denote by $P_{s,t}(n)$
the $n\times n$ submatrix 
of $P$ with coefficients ${i+j+s+t\choose i+s},\ 0\leq i,j<n$
and denote by $D_{s,t}(n)=\hbox{det}(P_{s,t}(n))$ its determinant.

{\bf Theorem 1.1.} {\sl We have
$$D_{s,t}(n)=\prod_{k=0}^{s-1}\frac{{n+k+t\choose t}}{{k+t\choose
t}}\quad ,s,t,n\geq 0\ .$$

In particular, the function $n\longmapsto D_{s,t}(n)$ is a polynomial of degree
$st$ in $n$.
}

This Theorem follows for instance from the formulas contained in section 5 of 
{\bf [GV]} (a beautiful paper studying mainly determinants of finite 
submatrices of the matrix
$T$ with coefficients $t_{i,j}={j\choose i}$). We give briefly a different proof
using the so-called \lq\lq condensation method'' (cf. for instance the survey 
paper {\bf [K1]}).

{\bf Proof of Theorem 1.1.} The definition ${a+b\choose a}=\frac{(a+b)!}{a!\quad b!}$
and a short computation show that Theorem 1.1 boils down to
$$\hbox{det}(A_k(n))=\prod_{i=0}^{n-1}i!\ (i+k)!$$
where $A_k(n)$ has coefficients
$a_{i,j}=(i+j+k)!$ for $0\leq i,j<n$ with $k=s+t$. The condensation identity
(cf. Proposition 10 in {\bf [K1]})
$$\hbox{det}(M)\ \hbox{det}(M_{1,n}^{1,n})=
\hbox{det}(M_1^1)\ \hbox{det}(M_n^n)-
\hbox{det}(M_1^n)\ \hbox{det}(M_n^1)$$
(with $M_{i_1,\dots,i_k}^{j_1,\dots,j_k}$ denoting the submatrix 
of the $n\times n$
matrix $M$ obtained by erasing lines $i_1,\dots,i_k$ and columns $j_1,\dots,j_k$)
allows a recursive (on $n$) computation of $\hbox{det}(A_k(n))$ 
establishing the result. \hfill QED

Theorem 1.1 has the following generalization. Let 
$$q(x,y)=\sum_{s,t=0}c_{s,t}{x\choose s}{y\choose t}\in{\bf Q}
[x,y]$$ be a polynomial in two variables $x,y$ and let $Q(n)$ be the matrix
with coefficients $q_{i,j}=q(i,j),\ 0\leq i,j<n$.

Elementary operations on rows and columns  show easily the following result.

{\bf Proposition 1.2.} {\sl One has for all $n$
$$\hbox{det}(P_{0,0}(n)+Q(n))=\hbox{det}(C_Q(n)+\hbox{Id}_n)$$
where $C_Q(n)$ has coefficients $c_{i,j}$ for $0\leq i,j<n$ and
where $\hbox{Id}_n$ denotes the identity matrix of order $n$.

In particular, the sequence of determinants
$$\hbox{det}(P_{0,0}(0))=1,\hbox{det}(P_{0,0}(1))=1+c_{0,0},\hbox{det}(P_{0,0}(2)),
\dots$$
becomes constant for $n\geq \mu$  where 
$\mu=\hbox{min}(\hbox{degree}_x(Q),\hbox{degree}_y(Q))$
with $\hbox{degree}_x(Q)$ (respectively $\hbox{degree}_x(Q)$)
denoting the degree of $Q$ with respect to $x$ (respectively $y$).}

In general, the function
$$n\longmapsto \hbox{det}(P_{s,t}(n)+Q(n))$$
seems to be polynomial of degree $\leq st$ in $n$ for $n$ huge enough.

Consider the symmetric matrix 
$G$ of order $k$ with coefficients $g_{i,j}=\sum_{s=0}^{n+k-1} {s\choose i}
{s\choose j}$ for $0\leq i,j<k$. Theorem 1.1 implies 
$\hbox{det}(G)=D_{k,k}(n)$ (with $D_{k,k}(n)$ given by the formula of
Theorem 1.1).

Let us also mention the following computation involving inverses of
binomial coefficients. Given three integers $s,t,n\geq 0$ let
$d_{s,t}(n)$ denote the determinant of the $n\times n$ matrix $M$
with coefficients
$$m_{i,j}={i+s+j+t\choose i+s}^{-1}\quad \hbox{for }0\leq i,j<n\ .$$

{\bf Theorem 1.3.} {\sl One has
$$d_{s,t}(n)=(-1)^{n\choose 2}\frac{1}{\prod_{k=0}^{n-1}
{2k+s+t\choose k+s}{2k-1+s+t\choose k}}\ .$$}

{\bf Sketch of proof.} For $0\leq k\in{\bf N}$ introduce the symmetric matrix
$A_k(n)$ of order $n$ with coefficients $a_{i,j}=\frac{1}{(i+j+k)!},\ 0\leq i,j<n$.
A small computation shows then that Theorem 1.6 is equivalent to the identity
$$\hbox{det}(A_k(n))=(-1)^{n\choose 2}\ \prod_{i=0}^{n-1}\frac{i!}{(n+k+i-1)!}$$
(with $k=s+t$) which can be proven recursively on $n$ by the
condensation method (cf. proof of Theorem 1.1).

Let us now consider the following variation of the Pascal triangle.
Recall that a complex matrix of rank $1$ and order $n\times n$ has
coefficients $\alpha_i\beta_j$ (for $0\leq i,j<n$) where 
$\alpha=(\alpha_0,\dots,\alpha_{n-1})$ and $\beta=(\beta_0,\dots,
\beta_{n-1})$ are two complex sequences, well defined up to 
$\lambda\alpha,\frac{1}{\lambda}\beta$ for $\lambda\in{\bf C}^*$.

Given two infinite sequences $\alpha=(\alpha_0,\alpha_1,\dots)$
and $\beta=(\beta_0,\beta_1,\dots)$ consider the $n\times n$ matrix $A(n)$
with coefficients $a_{i,j}=a_{i-1,j}+a_{i,j-1}+\alpha_i\beta_j$ for
$0\leq i,j<n$ (where we use the convention $a_{i,-1}=a_{-1,i}=0$ for
all $i$).

{\bf Proposition 1.4.} {\sl (i) The coefficient $a_{i,j}$ (for $0\leq i,j<n$) of
the matrix $A(n)$ is given by
$$a_{i,j}=\sum_{s=0}^i\sum_{t=0}^j\alpha_{i-s}\beta_{j-t}{s+t\choose s}\ .$$

\ \ (ii) The matrix $A(n)$ has determinant $(\alpha_0\beta_0)^n$.}

{\bf Proof.} Assertion (i) is elementary and left to the reader.

Assertion (ii) obviously holds if $\alpha_0=0$ or $\beta_0=0$.
We can hence suppose $\beta_0=1$. Proposition 1.1 and elementary 
operations on rows establish the result easily for arbitrary $\alpha$
and $\beta=(1,0,0,0,\dots)$. The case of an arbitrary sequence 
$\beta$ with $\beta_0=1$ is then reduced to the previous case using 
elementary operations on columns. \hfill QED

Another variation on theme of Pascal triangles is given by considering 
the $n\times n$ matrix $A(n)$ with coefficients $a_{i,0}=\rho^i,\ a_{0,i}
=\sigma^i,\
0\leq i<n$ and $a_{i,j}=a_{i-1,j}+a_{i,j-1}+x\ a_{i-1,j-1},\ 1\leq i,j<n$.
Setting $x=0,\ \rho=\sigma=1$ we get hence the matrix defined by
binomial coefficients considered above. One has then the following
result, due to C. Krattenthaler {\bf [K2]}, which we state without proof.

{\bf Theorem 1.5.} {\sl One has
$$\hbox{det}(A(n))=(1+x)^{n-1\choose 2}\ 
(x+\rho+\sigma-\rho\sigma)^{n-1}\ .$$}

Let now $B(n)$ be the symplectic (antisymmetric) $n\times n$ matrix defined by
$b_{i,i}=0,\ 0\leq i<n$, $b_{i,0}=-b_{0,i}=\rho^{i-1},\ 1\leq i<n$, $b_{i,j}=
b_{i-1,j}+b_{i,j-1}+x\ b_{i-1,j-1},\ 1\leq i,j<n$. 

The computation of the determinant of $B(2n)$ is again due to C. 
Krattenthaler {\bf [K2]}:
$$\hbox{det}(B(2n))=(1+x)^{2(n-1)^2}\ (\rho+x)^{2n-2}\ .$$

%\end{document}

\section{Generalized Pascal triangles}

Let $\alpha=(\alpha_0,\alpha_1,\dots)$ and $\beta=(\beta_0,\beta_1,\dots)$
be two sequences starting with a common first term 
$\gamma_0=\alpha_0=\beta_0$. 
Define a matrix $P_{\alpha,\beta}(n)$ of order 
$n$ with coefficients $p_{i,j}$ by setting 
$p_{i,0}=\alpha_i,
\ p_{0,i}=\beta_i$ for $0\leq i<n$ and $p_{i,j}=p_{i-1,j}+p_{i,j-1}$ 
for $1\leq i,j<n$. 

It is easy to see that the coefficient 
$p_{i,j}$ of
$P_{\alpha,\beta}(n)$ is also given by the formula
$$p_{i,j}=\gamma_0{i+j\choose i}+\left(\sum_{s=1}^i(\alpha_s-\alpha_{s-1})
{i-s+j\choose j}\right)+
\left(\sum_{t=1}^j(\beta_t-\beta_{t-1}){i+j-t\choose i}\right)\ .$$

We call the infinite \lq\lq matrix'' $P_{\alpha,\beta}(\infty)$ the {\it
generalized Pascal triangle} associated to $\alpha,\beta$.

We will mainly be interested in the sequence of determinants
$$(\hbox{det}(P_{\alpha,\beta}(1))=\gamma_0,
\hbox{det}(P_{\alpha,\beta}(2))=\gamma_0(\alpha_1+\beta_1)-\alpha_1\beta_1,
\dots,\hbox{det}(P_{\alpha,\beta}(n)),\dots)\ .$$

{\bf Example 2.1.} Take an arbitrary sequence $\alpha=(\alpha_0,\alpha_1,\dots)$
and let $\beta$ be the constant sequence 
$\beta=(\alpha_0,\alpha_0,\alpha_0,\dots)$. Proposition 1.5 implies 
$\hbox{det}(A_{(\alpha_0,\alpha_1,\dots),
(\alpha_0,\alpha_0,\dots)}(n))=\alpha_0^n$ (using perhaps the convention 
$0^0=1$). This yields an easy way of writing down matrices with
determinant $1$ by choosing a sequence $\alpha=(\alpha_0=1,\alpha_1,
\dots)$. The finite sequence $\alpha=(1,-2,5,11)$ for instance yields the
determinant $1$ matrix
$$\left(\begin{array}{cccc}
1&1&1&1\cr
-2&-1&0&1\cr
5&4&4&5\cr
11&15&19&24\end{array}\right)\ .$$

\section{Linear recursions}

This section is devoted to general Pascal triangles constructed from
sequences satisfying linear recursions. Conjecturally, the sequence
of determinants of such matrices satisfies then again a (generally much
longer) linear recursion. We prove this in the particular case where
the defining sequences are of order at most 2.

{\bf Definition.} A sequence $\sigma=(\sigma_0,\sigma_1,\sigma_2,\dots)$ 
satisfies a {\it linear recursion of order $d$} if there exist constants
$D_1,D_2,\dots,D_d$ such that
$$\sigma_n=\sum_{i=1}^d D_i\ \sigma_{n-i}\quad \hbox{ for all }n\geq d\ .$$
The polynomial 
$$z^d-\sum_{i=1}^dD_i\ z^{d-i}$$
is then called the {\it characteristic polynomial} of the linear recursion.

Let us first consider generalized Pascal triangles defined by linear 
recursion sequences of order at most $2$:

Given $\gamma_0,\alpha_1,\beta_1,A_1,A_2,B_1,B_2$ we set $\alpha_0=\beta_0=
\gamma_0$ and consider the square matrix $M(n)$ of order $n$ with entries
$$\begin{array}{l}
m_{i,0}=\alpha_i,\ 0\leq i<n\hbox{ where }\alpha_k=A_1\alpha_{k-1}+A_2
\alpha_{k-2},\ k\geq 2,\cr
m_{0,i}=\beta_i,\ 0\leq i<n\hbox{ where }\beta_k=B_1\beta_{k-1}+B_2
\beta_{k-2},\ k\geq 2,\cr
m_{i,j}=m_{i-1,j}+m_{i,j-1},\ 1\leq i,j<n\ .\end{array}$$
The matrix $M(3)$ for instance is hence given by
$$M(3)=\left(\begin{array}{ccc}
\gamma_0&\beta_1&B_1\beta_1+B_2\gamma_0\cr
\alpha_1&\alpha_1+\beta_1&\alpha_1+\beta_1+B_1\beta_1+B_2\gamma_0\cr
A_1\alpha_1+A_2\gamma_0&\alpha_1+\beta_1+A_1\alpha_1+A_2\gamma_0&
m_{3,3}\end{array}
\right)$$
where $m_{3,3}=
2\alpha_1+2\beta_1+A_1\alpha_1+B_1\beta_1+A_2\gamma_0+B_2\gamma_0$.

We have hence $M(n)=P_{\alpha,\beta}(n)$ where $P_{\alpha,\beta}$
is the generalized Pascal triangle introduced in the previous section.

We set $d(n)=\hbox{det}(M(n))$ for $n\geq 1$ and introduce the constants
$$\begin{array}{l}
D_1=-(A_1\beta_1+B_1\alpha_1-2(\alpha_1+\beta_1)+\gamma_0\left(
A_1B_2+A_2B_1-(A_2+B_2)+A_2B_2\right))  \cr
D_2=-\left(A_2\gamma_0+\alpha_1+(1-A_1-A_2)\beta_1\right)
\left(B_2\gamma_0+\beta_1+(1-B_1-B_2)\alpha_1\right) \ .\end{array}$$

{\bf Theorem 3.1.} {\sl The sequence $d(n),\ n\geq 1$ defined as above
satisfies the following equalities
$$\begin{array}{l}
d(1)=\gamma_0 \ ,\cr
d(2)=\gamma_0(\alpha_1+\beta_1)-\alpha_1\beta_1\ ,\cr
d(n)=D_1\ d(n-1)+D_2\ d(n-2)\hbox{ for all }n\geq 3\ .\end{array}$$}

Theorem 3.1 will be proven below.

{\bf Example 3.2} (a) 
The sequence $(\hbox{det}(P_{\alpha,\beta}(n)))_{n=1,2,\dots}$
of determinants associated to two geometric sequences
$$\begin{array}{l}
\alpha=(1,A,A^2,\dots, \alpha_k=A^k,\dots)\cr
\beta=(1,B,B^2,\dots, \beta_k=B^k,\dots)
\end{array}$$
is given by
$$\hbox{det}(P_{\alpha,\beta}(n))=(A+B-AB)^{n-1}\ .$$

Let $\alpha=(\alpha_0,\alpha_1,\dots)$ and $\beta=(\beta_0,\beta_1,\dots)$ be
two sequences satisfying $\alpha_0=\beta_0=\gamma_0$ and linear recursions
$$\begin{array}{l}
\alpha_n=\sum_{i=1}^a A_i\alpha_{n-i}\hbox{ for }n\geq a\ ,\cr
\beta_n=\sum_{i=1}^b B_i\beta_{n-i}\hbox{ for }n\geq b\end{array}$$
of order $a$ and $b$.

Theorem 3.1 and computations suggest that the following might be true.

{\bf Conjecture 3.3.} {\sl If two sequences $\alpha=(\alpha_0,\alpha_1,\dots),
\beta=(\beta_0,\beta_1,\dots)$ satisfy both linear recurrence relations then 
there exist a natural integer $d\in{\bf N}$ and
constants $D_1,\dots,D_d$ (depending on $\alpha,\beta$) such that
$$\hbox{det}(P_{\alpha,\beta}(n))=\sum_{i=1}^d D_i\
\hbox{det}(P_{\alpha,\beta}(n-i))\quad \hbox{ for all }n>d\ .$$}

{\bf Remarks 3.4.} (i) Generically, (ie. for $\alpha$ and $\beta$
two generic sequences of order $a$ and $b$ such that $\alpha_0=\beta_0$)
the integer $d$ of Conjecture 3.3 seems to be given by $d={a+b-2\choose a-1}$.

\ \ (ii) Generically, the coefficient $D_i$ seems to be 
a homogeneous form 
(with polynomial coefficients in $A_1,\dots,A_a,B_1,B_b$) of degree $i$
in $\gamma_0,\alpha_1,\alpha_{a-1},\beta_1,\dots\beta_{b-1}$. For 
non-generic pairs of sequences (try $\beta=-\alpha$ with 
$\alpha=(0,\alpha_1,\dots)$ satisfying a linear recursion of order 3)
the coefficients $D_i$ may be rational fractions in the variables.

\ \ (iii) If $a=b>1$ and the recursive sequences $\alpha,\beta$ are
generic, then the coefficients $D_0=-1,D_1,\dots,D_d$ 
of the linear recursion in Conjecture 3.3 seem to have the symmetry
$$D_{d-i}=q^{(d-2i)/2}D_i$$
where
$q$ is a quadratic form in $\gamma_0,\alpha_i,\beta_i$ factorizing
into a product of two linear forms which are symmetric under the exchange
of parameters
$\alpha_i$ with $\beta_i$ and $A_i$ with $B_i$ (this corresponds to
transposing $P_{\alpha,\beta}$).

Theorem 3.1 shows that the generic quadratic form $q_2$ working for 
$a=b=2$ is given by
$$q_2=\left(A_2\gamma_0+\alpha_1+(1-A_1-A_2)\beta_1\right)
\left(B_2\gamma_0+\beta_1+(1-B_1-B_2)\alpha_1\right)\ .$$
The generic quadratic form $q_3$ working for $a=b=3$ seems to be
$$\begin{array}{lll}
\displaystyle 
q_3&\displaystyle=&\displaystyle
\left(A_3\gamma_0+\alpha_1+\alpha_2+(1-A_1-A_2-A_3)\beta_1\right)\cr
&&\displaystyle \qquad
\left(B_3\gamma_0+\beta_1+\beta_2+(1-B_1-B_2-B_3)\alpha_1\right)\ .
\end{array}$$

{\bf Example 3.5.} Consider the $3-$periodic sequence 
$\alpha=(\alpha_0,\alpha_1,\alpha_2,\dots,\alpha_k=\alpha_{k-3},\dots)$. 
The sequence $d(n)=\hbox{det}(P_{\alpha,\alpha}(n))$ seems then to satisfy 
the recursion relation
$$\begin{array}{l} d(n)=D_1\ d(n-1)+D_2\ d(n-2)-(\alpha_0+\alpha_1+\alpha_2)
D_2\ d(n-3)\cr
\quad -(\alpha_0+\alpha_1+\alpha_2)^3 D_1\ d(n-4)+(\alpha_0+\alpha_1+\alpha_2)^5\
d(n-5)\end{array}$$
where
$$\begin{array}{l}
D_1=11\alpha_1+5\alpha_2\cr
D_2=-(3\alpha_0^2+37\alpha_1^2+3\alpha_2^2+15\alpha_0\alpha_1+
5\alpha_0\alpha_2+24\alpha_1\alpha_2)\ .\end{array}$$

In the general case
$$\begin{array}{l}
\alpha=(\alpha_0=\gamma_0,\alpha_1,\alpha_2,\dots,\alpha_k=\alpha_{k-3},\dots)\cr
\beta=(\beta_0=\gamma_0,\beta_1,\beta_2,\dots,\beta_k=\beta_{k-3},\dots)
\end{array}$$
of two $3-$ periodic sequences (starting with a common value $\gamma_0$) one 
seems to have
$$d(n)=D_1d(n-1)+D_2d(n-2)+D_3d(n-3)+q D_2d(n-4)+q^2 D_1d(n-5)-q^3
d(n-6)$$
where 
$$q=(\gamma_0+\alpha_1+\alpha_2)(\gamma_0+\beta_1+\beta_2)$$
$$D_1=\gamma_0+6(\alpha_1+\beta_1)+3(\alpha_2+\beta_2)$$
$$\begin{array}{l}
D_2=-\big(3\gamma_0^2+12(\alpha_1^2+\beta_1^2)+13
\gamma_0(\alpha_1+\beta_1)+5\gamma_0(\alpha_2+\beta_2)\cr
\qquad +9(\alpha_1\alpha_2+\beta_1\beta_2)+11(\alpha_1\beta_2+\alpha_2\beta_1)
+24\alpha_1\beta_1+8\alpha_2\beta_2\big)\end{array}$$
$$\begin{array}{l}
D_3=6\gamma_0^3+\gamma_0^2(18(\alpha_1+\beta_1)+8(\alpha_2+\beta_2))
+\gamma_0(25(\alpha_1^2+\beta_1^2)+3(\alpha_2^2+\beta_2^2)\cr
\quad +
18(\alpha_1\alpha_2+\beta_1\beta_2)+54\alpha_1\beta_1+26(\alpha_1\beta_2+\alpha_2\beta_1)+10\alpha_2\beta_2)\cr
\quad +9(\alpha_1^3+\beta_1^3)+9(\alpha_1^2\alpha_2+\beta_1^2\beta_2)
+28(\alpha_1^2\beta_1+\alpha_1\beta_1^2)+22(\alpha_1^2\beta_2+\alpha_2\beta_1^2)
\cr
\quad
+3(\alpha_1\beta_2^2+\alpha_2^2\beta_1)+3(\alpha_2^2\beta_2+\alpha_2\beta_2^2)
+30(\alpha_1\alpha_2\beta_1+\alpha_1\beta_1\beta_2)+
24(\alpha_1\alpha_2\beta_2+\alpha_2\beta_1\beta_2)
\end{array}$$

Let us briefly explain how Conjecture 3.3 can
be tested on a given pair $\alpha,\beta$ of linear recurrence sequences.

First Step. Guess $d$.

Second step. Compute at least $2d+1$ terms of the sequence
$$w_1=\hbox{det}(P_{\alpha,\beta}(1)),w_2=\hbox{det}(P_{\alpha,\beta}(2)),
\dots\ .$$

Third step. Check that the so-called {\it Hankel matrix}
$$H_{d+1}(w)=\left(\begin{array}{ccccc}
w_1&w_2&w_3&\dots&w_{d+1}\cr
w_2&w_3&w_4&\dots&w_{d+2}\cr
\vdots\cr
w_{d+1}&w_{d+2}&w_{d+3}&\dots&w_{2d+1}\end{array}\right)$$
of order $d+1$ has zero determinant (otherwise try again with a higher value
for $d$) and choose a vector of the form
$$L=(D_d,D_{d-1},D_{d-2},\dots,D_2,D_1,-1)$$
in its kernel. One has then by definition
$$\hbox{det}(P_{\alpha,\beta}(n))=\sum_{i=1}^d D_i\ 
\hbox{det}(P_{\alpha,\beta}(n-i))$$
for $d+1\leq n\leq 2d+1$.

Finally, check (perhaps) the above recursion for a few more values of $n>2d+1$.

\subsection{Proof of Theorem 3.1.}

The assertions concerning $d(1)$ and $d(2)$ are obvious.
One checks (using for instance a symbolic computation program on a computer) 
that the recursion relation holds for $d(3),\ d(4)$ and $d(5)$.

Introduce now the lower and upper triangular square matrices
$$T_A=\left(\begin{array}{ccccccccc}
1&0&0&0&\dots\cr
-A_1&1&0&0&\dots\cr
-A_2&-A_1&1&0&\dots\cr
0&-A_2&-A_1&1&\dots\cr
\vdots&&\ddots\end{array}\right)$$
$$T_B=\left(\begin{array}{ccccccccc}
1&-B_1&-B_2&0&0&\dots\cr
0&1&-B_1&-B_2&0&\dots\cr
0&0&1&-B_1&-B_2&\dots\cr
&\vdots&&\ddots\end{array}\right)$$
of order $n$ 
and set $\tilde M=T_A\ M\ T_B$. The entries $\tilde m_{i,j},\ 0\leq i,j<n$ of
$\tilde M$ satisfy $\tilde m_{i,j}=\tilde m_{i-1,j}+\tilde m_{i,j-1},\ (i,j)\not=(2,2)$ for 
$2\leq i,j<n$. One has 
$$\tilde M=\left(\begin{array} {cccccc}
\gamma_0&\beta_1-B_1\gamma_0&0&0&0&\dots \cr
\alpha_1-A_1\gamma_0&\tilde m_{1,1}&\tilde m_{1,2}&\tilde m_{1,3}&
\tilde m_{1,4}&\dots\cr 
0&\tilde m_{2,1}&\tilde m_{2,2}&\tilde m_{2,3}&
\tilde m_{2,4}&\dots\cr
\vdots&&\vdots\end{array}\right)$$
where 
$$\begin{array}{l}
\tilde m_{1,1}=\alpha_1+\beta_1-A_1\beta_1-B_1\alpha_1+A_1B_1\gamma_0 \cr
\tilde m_{1,2}=\tilde m_{1,3}=\tilde m_{1,4}=\dots=(
1-B_1-B_2)\alpha_1+\beta_1+B_2\gamma_0\cr
\tilde m_{2,1}=\tilde m_{3,1}=\tilde m_{4,1}=\dots=(
1-A_1-A_2)\beta_1+\alpha_1+A_2\gamma_0\cr
\tilde m_{2,2}=(2-B_1)\alpha_1+(2-A_1)\beta_1+(A_2+B_2-A_1B_2-A_2B_1-A_2B_2)
\gamma_0
\end{array}$$
Developping the determinant $d(n)=\hbox{det}(\tilde M)$ along the
second row of $\tilde M$ one obtains
$$d(n)=(\gamma_0(\alpha_1+\beta_1)-\alpha_1\beta_1)
\overline d(n-2)+\gamma_0\ \hbox{det}(P(n-1))$$
where $\overline d(n-2)=\hbox{det}(\overline M(n-2))$ with coefficients
$\overline m_{i,j}=\tilde m_{i+2,j+2}$ for $0\leq i,j<n-2$ (ie. $\overline
M(n-2)$ is the principal submatrix of $\tilde M$ obtained by erasing
the first two rows and columns of $\tilde M$) and where $P(n-1)$
is the square matrix of order $(n-1)$ with entries $p_{0,0}=0$ and
$p_{i,j}=\tilde m_{i+1,j+1}$ for $0\leq i,j<n-1, (i,j)\not=(0,0)$.

The matrix $\overline M(m)$ ($m\leq n-2$) is a generalized Pascal triangle 
associated to the linear recursion sequences 
$\overline \alpha=(\overline \alpha_0,\overline \alpha_1,\dots)$
and $\overline \beta=(\overline \beta_0,\overline \beta_1,\dots)$
of order $2$ defined by
$$\begin{array}{l}
\overline \alpha_0=\overline \beta_0=
(2-B_1)\alpha_1+(2-A_1)\beta_1+(A_2+B_2-A_1B_2-A_2B_1-A_2B_2)
\gamma_0
 ,\cr
\overline \alpha_1=(3-B_1)\alpha_1+(3-2A_1-A_2)\beta_1+
(2A_2+B_2-A_1B_2-A_2B_1-A_2B_2)\gamma_0  \cr
\overline \beta_1=(3-2B_1-B_2)\alpha_1+(3-A_1)\beta_1+
(A_2+2B_2-A_1B_2-A_2B_1-A_2B_2)\gamma_0 \cr
\overline \alpha_n=2\overline \alpha_{n-1}-\overline \alpha_{n-2}\quad
\hbox{for }n\geq 2\ ,\cr
\overline \beta_n=2\overline \beta_{n-1}-\overline \beta_{n-2}\quad
\hbox{for }n\geq 2\ .\end{array}$$
Induction on $n$ and a computation (with
$\overline A_1=\overline B_1=2,\ \overline A_2=\overline B_2=-1$) shows
the equality
$$\overline d(m)=D_1 \ \overline d(m-1)+D_2\ \overline d(m-2)$$
for $3\leq m<n$ where $D_1$ and $D_2$ are as in the Theorem.

Introducing the $(n-1)\times (n-1)$ lower triangular square matrix
$$T_P=\left(\begin{array}{ccccccccc}
1&0&0&\dots\cr
-1&1&0&\dots\cr
0&-1&1&\dots\cr
&\vdots&&\ddots\end{array}\right)$$
we get $\tilde P=T_P\ P(n-1)\ T_P^t$ with coefficients $\tilde p_{i,j},\ 
0\leq i,j<n-1$ given by 
$$\begin{array}{l}
\tilde p_{0,0}=\tilde p_{i,0}=\tilde p_{0,i}=0 \hbox{ for }2\leq i<n-1,\cr
\tilde p_{0,1}=(1-B_1-B_2)\alpha_1+\beta_1+B_2\gamma_0 \cr
\tilde p_{1,0}=(1-A_1-A_2)\beta_1+\alpha_1+A_2\gamma_0 \cr
\tilde p_{1,i}=\tilde p_{0,1}  \hbox{ for }2\leq i<n-1,\cr
\tilde p_{i,1}=\tilde p_{1,0} \hbox{ for }2\leq i<n-1,\cr
\tilde p_{2,2}=(2-B_1)\alpha_1+(2-A_1)\beta_1+(A_2+B_2-A_1B_2-A_2B_1-A_2B_2)
\gamma_0 \cr
\tilde p_{i,j}=\tilde p_{i-1,j}+\tilde p_{i,j-1}\ ,2\leq i,j<n-1,\ (i,j)\not=
(2,2)\ .\end{array}$$

Let $\overline P(n-3)$ denote the square matrix of order $(n-3)$ with 
coefficients $\overline p_{i,j}=\tilde p_{i+2,j+2}\ ,0\leq i,j<n-3$ (ie.
$\overline P(n-3)$ is obtained by erasing the first two rows and columns 
of $\tilde P(n-1)$). One checks the equality
$$\overline P(n-3)=\overline M(n-3)$$
where $\overline M(n-3)$ is defined as above. This implies the identity
$$\begin{array}{l}
d(n)=(\gamma_0(\alpha_1+\beta_1)-\alpha_1\beta_1)\overline d(n-2)\cr
\quad -\gamma_0((1-B_1-B_2)\alpha_1+\beta_1+B_2\gamma_0)((1-A_1-A_2)\beta_1
+\alpha_1+A_2\gamma_0)\overline d(n-3)\ .\end{array}$$

Using the recursion relation 
$\overline d(m)=D_1 \ \overline d(m-1)+D_2\ \overline d(m-2)$
(which holds by induction for $3\leq m<n$) we can hence express $d(n)$ as a 
linear function
(with polynomial coefficients in $\gamma_0,\alpha_1,\beta_1,A_1,A_2,
B_1,B_2$) of $\overline d(m-4)$ and $\overline d(m-5)$.

Comparison of this with the linear expression in 
$\overline d(m-4)$ and $\overline d(m-5)$ obtained similarly from
$D_1\ d(n-1)+D_2\ d(n-2)$ finishes the proof.

\section{Symmetric matrices}

Take an arbitrary sequence $\alpha=(\alpha_0,\alpha_1,\dots)$. The generalized
Pascal triangle associated to the pair of identical sequences $\alpha,\alpha$
is the {\it generalized symmetric Pascal triangle} associated to 
$\alpha$ and yields symmetric matrices $P_{\alpha,\alpha}(n)$
by considering principal submatrices consisting
of the first $n$ rows and columns of $P_{\alpha,\alpha}$. 

The main example is of course the classical Pascal triangle obtained
from the constant sequence $\alpha=(1,1,1,\dots)$. Other sequences
satisfying linear recursions like for instance the Fibonacci sequence
$$(0,1,1,2,3,5,8,\dots)$$
and shifts of it yield also nice examples. 

Conjecture 3.3 should of course also hold for symmetric matrices 
obtained by considering the generalised symmetric Pascal triangle 
associated to a sequence satisfying a linear recurrence relation.

The generic order $d_s(a)$ (where $a$ denotes the order of the defining 
linear sequence) of the linear recursion satisfied by
$\hbox{det}(P_{\alpha,\alpha}(n))$
seems however usually to be smaller than in the generic
non-symmetric case. Examples yield the following first values
$$\begin{array}{lccccccc}
a=\quad&1&2&3&4&5&6\cr
d_s(a)=\quad&1&2&5&14&41&122
\end{array}$$
and suggest that perhaps $d_a=(3^{a-1}+1)/2$.

The coefficients $D_i$ seem still to be polynomial in $\alpha_i$ and
$A_i$.

The symmetry relation has also an analogue (in the generic case)
which is moreover
somewhat simpler in the sense that it is given by a linear form $\rho$ 
(in $\alpha_0,\dots,\alpha_{a-1}$) and we seem to have
$$D_{d_s-i}=\rho^{d_s-2i}\ D_i$$
(where $D_0=-1$).

{\bf Example 4.1.} If a sequence 
$$\alpha=(\alpha_0,\alpha_1,\alpha_2,\dots,\alpha_k=A_1\alpha_{k-1}+A_2\alpha_{k-2}
+A_3\alpha_{k-3},\dots)$$
satisfies a linear recursion relation of order $3$, then the sequence
$d(n)=\hbox{det}(P_{\alpha,\alpha}(n))$ (the
matrix
$P_{\alpha,\alpha}(n)$ has coefficients $p_{0,i}=p_{i,0}=\alpha_i,\ 0\leq i<n$ and
$p_{i,j}=p_{i-1,j}+p_{i,j-1}$ for $1\leq i,j<n$)
of the associated determinants seems to satisfy
$$\begin{array}{l}
d(1)=\alpha_0\ ,\cr
d(2)=2\alpha_0\alpha_1-\alpha_1^2\ , \cr
d(3)=(2\alpha_1-\alpha_2)\left(\alpha_0(2\alpha_1+\alpha_2)-2\alpha_1^2\right)
)\ ,\cr
d(n)=D_1\ d(n-1)+D_2\ d(n-2)+\rho D_2\ d(n-3)+\rho^3\ D_1\ d(n-4)-\rho^5
\ d(n-5)\end{array}$$
where
$$\begin{array}{l}
\rho=-A_3\alpha_0+(-2+2A_1+A_2+A_3)\alpha_1-\alpha_2 \ ,\cr
D_1=A_3(1-2A_1-2A_2-A_3)\alpha_0+(10-10A_1-A_2+A_3+4A_1^2+2A_1A_2)\alpha_1\cr
\qquad +(5-4A_1-2A_2)\alpha_2  \ ,\cr
D_2=c_{0,0}\alpha_0^2+c_{1,1}\alpha_1^2+c_{2,2}\alpha_2^2+c_{0,1}\alpha_0\alpha_1+
c_{0,2}\alpha_0\alpha_2+c_{1,2}\alpha_1\alpha_2\end{array}$$
with
$$\begin{array}{l}
c_{0,0}=-A_3^2(2-2A_1+2A_2+A_3+A_1^2)\ ,\cr
c_{1,1}=-40+80A_1+16A_2+4A_3\cr
\qquad-64A_1^2-2A_2^2-A_3^2-28A_1A_2-20A_1A_3-2A_2A_3\cr
\qquad +2A_1(2A_1+A_2+A_3)(6A_1+A_2+A_3)-A_1^2(2A_1+A_2+A_3)^2\ ,\cr
c_{2,2}=-10+12A_1+6A_2+8A_3-(2A_1+A_2+A_3)^2\ ,\cr
c_{0,1}=-A_3(16-28A_1+16A_1^2-2A_2^2-A_3^2+2A_1A_3-3A_2A_3\cr
\qquad -2A_1^2(2A_1+A_2+A_3))\ ,\cr
c_{0,2}=-A_3(8-10A_1-3A_3+2A_1(2A_1+A_2+A_3))\ ,\cr
c_{1,2}=2\left(-20+32A_1+10A_2+9A_3-18A_1^2-A_2^2-A_3^2\right.\cr
\qquad \left.-11A_1A_2-12A_1A_3-2A_2A_3+A_1(2A_1+A_2+A_3)^2\right)\end{array}$$

We conclude this section by mentionning the following 
more exotic example:

{\bf Example 4.2.} (Central binomial coefficients)
Consider the sequence 
$$\alpha=(1,2,6,20,70,252,924,3432,12870,48620,184756,\dots,\alpha_k={2k\choose
k},\dots)$$
of central binomial coefficients. For $1\leq n\leq 36$ the values of 
$\hbox{det}(P_{\alpha,\alpha}(n))$ are zero except if $n\equiv 1,3\pmod 6$ and
for $n\equiv 1,3\pmod 6$ the values of $\hbox{det}(P_{\alpha,\alpha}(n))$ have the
following intriguing factorisations:
$$\begin{array}{llcll}
n:&\hbox{det}(P_{\alpha,\alpha}(n))&\quad&n:&\hbox{det}(P_{\alpha,\alpha}(n))\cr
1:&1& &3:&-2^2   \cr
7:&-2^6& &9:&2^8\ 3^2  \cr
13:&2^{16}\ 3^6& &15:&-2^{26}\ 3^6  \cr
19:&-2^{30}\ 3^2\ 103^4& &21:&2^{24}\ 103^4\ 4229^2  \cr
25:&2^{24}\ 3^4\ 31^4\ 431^4\ 4229^2& &27:&-2^{26}\ 3^6\ 31^4\ 59^2\ 431^4\
11701^2  \cr
31:&-2^{30}\ 3^{10}\ 59^2\ 11701^2\ p^4& &33:&2^{32}\ 3^{12}\ 11^2\ 
2017^2\ 28349^2\ p^4
\end{array}$$
where $p=4893589$.

The matrix $P_{\alpha,\alpha}(n)$ seems to have rank $n$ if $n\equiv 1,3\pmod 6$,
rank $n-1$ if $n\equiv 0 \pmod 2$ and rank $n-2$ if $n\equiv 5\pmod 6$.

\section{Symplectic matrices}

Given an arbitrary sequence $\alpha=(\alpha_0,\alpha_1,\dots)$
with $\alpha_0=0$, the
matrices $P_{\alpha,-\alpha}(n)$ are symplectic (antilinear).

Determinants of integral symplectic matrices are squares of
integers and are zero in odd dimensions. We restrict hence ourself to even 
dimensions and consider sometimes also the (positive) square-roots of the
determinants. Even if Conjecture 3.3 holds there is of course no reason 
that the square roots of the determinants satisfy a linear recursion.

The conjectural recurrence relation for symplectic matrices has
the form
$$\hbox{det}(P_{\alpha,-\alpha}(2n))=\sum_{i=1}^{d(\alpha)}D_i
\hbox{det}(P_{\alpha,-\alpha}(2n-2i))\ .$$
However the coefficients $D_1,\dots ,D_{d({\alpha})}$ seem no longer to
be polynomial but rational for generic $\alpha$. Moreover, the nice
symmetry properties of the coefficients $D_i$
present in the other cases seem to have disappeared too.

{\bf Proposition 5.1.} {\sl 
(i) The symplectic matrices $P_{\alpha,-\alpha}(2n)$
associated to the sequence $\alpha=(0,1,1,1,1,\dots)$ have determinant
$1$ for every natural integer $n$.

\ \  (ii) The symplectic matrices $P_{\alpha,-\alpha}(2n)$
associated to the sequence $\alpha=(0,1,2,3,4,5,\dots)$ have determinant
$1$ for every natural integer $n$.}

Both assertions follow of course from Theorem 3.1. We will
however reprove them independently. 

{\bf Proof.} Consider the generalized Pascal triangle $$P=P_{(1,1,1,
\dots,1,\dots),(1,-1,-1,\dots,-1,\dots)}(\infty)$$
$$=\left(\begin{array}{cccccccc}
1&-1&-1&-1&-1&-1&-1&\dots\cr
1&0&-1&-2&-3&-4&-5&\dots\cr
1&1&0&-2&-5&-9&-14&\dots\cr
1&2&2&0&-5&-14&-28&\dots\cr
\vdots&\end{array}\right)\ .$$
The matrices $P(m)$ given by retaining only the first $m$ rows and columns of
$P(\infty)$ are all of determinant $1$ (compare the transposed
matrix $P(m)^t$ with Example 2.1).

Expanding the determinant along the first row one gets
$$1=\hbox{det}(P(m))$$
$$=\hbox{det}(P_{(0,1,1,1,\dots),-(0,1,1,\dots)}(m))+
\hbox{det}(P_{(0,1,2,3,4,\dots),-(0,1,2,3,\dots)}(m-1))\ .$$
The fact that symplectic matrices of odd order have zero determinant proves
now assertion (i) for even $m$ and assertion (ii) for odd $m$.\hfill QED

{\bf Remark 5.2.} The coefficients $p_{i,j}$ of the infinite symplectic matrix 
$$P_{(0,1,1,1,\dots),-(0,1,1,1,\dots)}(\infty)$$
have many interesting properties: One can for instance easily check
that
$$p_{i,j}={i+j-1\choose j}-{i+j-1\choose j-1}$$
(with the correct definition for 
${k\choose -1}$ given by ${-1\choose -1}={-1\choose 0}=1$ and
${k\choose -1}=0$ for $k=0,1,2,3,\dots$).
These numbers are the orders of the irreducible matrix algebras in the
Temperley-Lieb algebras (see for instance chapter 2.8, pages 86-101,
in {\bf [GHJ]}). 

There are other matrices constructed
using the numbers ${i+j-1 \choose j}-{i+j-1\choose j-1}$ whose determinants
have interesting properties: Let $A_k(n)$ and $B_k(n)$ be the 
$n\times n$ matrices with entries
$$a_{i,j}={2i+2j+k\choose i}-{2i+2j+k\choose i-1}\hbox{ and }b_{i,j}=
{2i+2j+k\choose i+1}-{2i+2j+k\choose i}$$
for $0\leq i,j<n$ and $k$ a fixed integer. C. Krattenthaler {\bf [K2]} has
shown that one has 
$$\hbox{det}(A_k(n))=2^{n\choose 2}$$
and 
$$\hbox{det}(B_k(n))=2^{n\choose 2}\frac{\prod_{i=0}^{n-1}(k+2i-1)}{n!}\ .$$

Principal minors of $P_{(0,1,1,1,\dots),-(0,1,1,1,\dots)}(\infty)$
associated to submatrices
consisting of $2n$ consecutive rows and columns and starting at rows 
and columns of index $k=0,1,2,\dots$ have interesting properties
as given by the following result (cited without proof)
which is an easy corollary of the computation of 
$\hbox{det}((a+j-i)\Gamma(b+i+j))\ ,\ 0\leq i,j<n$
by Mehta and Wang (cf. {\bf [MW]}).

{\bf Theorem 5.3.} {\sl Denote by $T_k(2n)$ the $2n\times 2n$ symplectic
matrix with coefficients 
$$t_{i,j}={2k+i+j-1\choose k+j}-{2k+i+j-1\choose 
k+j-1}=\frac{(i-j)(2k+i+j-1)!}{(k+i)!(k+j)!}$$ for $0\leq i,j<2n$. 
One has 
$$\sqrt{\hbox{det}(T_k(2n))}=\prod_{t=1}^{k-1}\frac{{2n+2t\choose t}}{{2t\choose t}}\ ,n=0,1,2,\dots\ .$$}

% $$\sqrt{\hbox{det}(T_k(2n))}=\prod_{t=0}^{k-2}\prod_{i=0}^t\frac{2n+2+t+i}{2+t+i}\ ,n=0,1,2,\dots\ .$$}

The first polynomials $$D_k(n)=\prod_{t=1}^{k-1}\frac{{2n+2t\choose t}}{{2t\choose t}}=\sqrt{\hbox{det}(T_k(2n))}$$ are given by
$$\begin{array}{lcl}
D_0(n)&=&1\cr
D_1(n)&=&1\cr
D_2(n)&=&(n+1)\cr
D_3(n)&=&(2n+3)(n+1)(n+2)/6\cr
D_4(n)&=&(2n+5)(2n+3)(n+3)(n+2)^2(n+1)/180\cr
\end{array}$$
 
The sequences $(D_k(n))_{k=0,1,2,\dots}$ (for fixed $n$) seem also to 
be of interest since they have appeared elsewhere. They start as follows:
$$\begin{array}{lcl}
(D_0(0),D_1(0),D_2(0),\dots)&=&(1,1,1,\dots)\cr
(D_0(1),D_1(1),D_2(1),\dots)&=&(1,1,2,5,14,\dots)\hbox{ (Catalan numbers)}
\cr
(D_0(2),D_1(2),D_2(2),\dots)&=&(1,1,3,14,84,\dots)\hbox{ (cf {\bf A}005700 in 
{\bf [IS]})}\cr
(D_0(3),D_1(3),D_2(3),\dots)&=&(1,1,4,30,330,\dots)\hbox{ (cf {\bf A}006149 in 
{\bf [IS]})}\cr
(D_0(4),D_1(4),D_2(4),\dots)&=&(1,1,5,55,1001,\dots)\hbox{ (cf {\bf A}006150 in 
{\bf [IS]})}\cr
(D_0(5),D_1(5),D_2(5),\dots)&=&(1,1,6,91,2548,\dots)\hbox{ (cf {\bf A}006151 in 
{\bf [IS]})}
\end{array}$$

Geometric sequences provide other nice special cases of Theorem 3.1.

{\bf Example 5.4.} (i) The sequence $\alpha=(0,1,A,A^2,A^3,\dots)$ (for $A>0)$
yields $\hbox{det}(P_{\alpha,-\alpha}(2n))=A^{2(n-1)}$.

\ \ (ii) The slightly more general example $\alpha=(0,1,A+B,\dots,
\alpha_k=\frac{A^k-B^k}{A-B},\dots)$ yields
$\hbox{det}(P_{\alpha,-\alpha}(2n))=(A-AB+B)^{2(n-1)}$.

Finally, we would like to mention the following exotic example.

{\bf Example 5.5.} The sequences 
$$\begin{array}{l}
\alpha_C=(0,1,1,2,5,14,42,\dots)\cr
\alpha_B=(0,1,2,6,20,70,\dots)\end{array}$$
related to Catalan numbers and central binomial coefficients
yield the sequences $r_C(n)=\sqrt{\hbox{det}(P_{\alpha_C,
-\alpha_C}(2n))}$ and 
$r_B(n)=\sqrt{\hbox{det}(P_{\alpha_B,-\alpha_B}(2n))}$:
$$\begin{array}{l}
\begin{array}{lcccccccccccccccccccccccccccccc}
n=&1&2&3&4&5&6&7\cr
r_C(n)=&1&2&6&31&286&4600&130664\cr
r_B(n)=&1&2\cdot 2&6\cdot 2^2&31\cdot 2^3&286\cdot 2^4&4600\cdot 2^5&
130664\cdot 2^6\end{array}\cr
\cr
\begin{array}{lcccccccccccccccccccccccccccccc}
n=&8&9&10&\dots\cr
r_C(n)=&6619840&591478944&93683332808&\dots\cr
r_B(n)=&6619840\cdot 2^7&591478944\cdot 2^8&93683332808\cdot 2^9&\dots\cr
\end{array}\end{array}$$

suggesting the conjecture $r_B(n)=2^{n-1} r_C(n)$ for $n\geq 1$.

\subsection{The even symplectic construction and the even
symplectic unimodular tree}

Given an arbitrary sequence
$\beta=(\beta_0,\beta_1,\dots)$ we consider the sequence 
$\alpha=(0,\beta_0,0,\beta_1,0,\beta_2,\dots)$ defined by $\alpha_{2n}=0$
and $\alpha_{2n+1}=\beta_n$. We call this way of constructing a
symplectic matrix $P_{\alpha,-\alpha}(2n)$ 
out of a sequence $\beta=(\beta_0,\beta_1,\dots)$ the
{\it even symplectic construction (of Pascal triangles)}.

{\bf Example 5.1.1.} The symplectic matrix of order $6$ associated to 
the the sequence $\beta=(1,1,-1,\dots)$ by the even symplectic construction
is the following determinant $1$ matrix
$$\left(\begin{array}{cccccc}
0&-1&0&-1&0&1\cr
1&0&0&-1&-1&0\cr
0&0&0&-1&-2&-2\cr
1&1&1&0&-2&-4\cr
0&1&2&2&0&-4\cr
-1&0&2&4&4&0\end{array}\right)\ .$$

By elementary operations on rows and columns it is easy to check the
identity
$$\begin{array}{l}
\hbox{det}(P_{(0,\beta_0,0,\beta_1,0,\beta_2,\dots),-(0,\beta_0,0,\beta_1,
\dots)}(2n))\cr
\qquad =
\hbox{det}(P_{(0,\beta_0,\beta_0,\beta_1,\beta_1,\beta_2,\beta_2,\dots),
-(0,\beta_0,\beta_0,\beta_1,\beta_1,\dots)}(2n))\end{array}$$
for all $n$ and $\beta=(\beta_0,\beta_1,\dots)$.

The main feature of the even symplectic construction is perhaps given by the
following result.

{\bf Theorem 5.1.2.} {\sl (i)
Let $(\beta_0,\beta_1,\dots,\beta_{n-1})$ be a sequence
of integers such that 
$$\hbox{det}(P_{(0,\beta_0,0,\beta_1,\dots,0,\beta_{n-1}),-
(0,\beta_0,0,\beta_1,\dots,0,\beta_{n-1})}(2n))=1\ .$$
Then there exists a unique even integer $\tilde \beta_n$ such that
$$\begin{array}{l}
\displaystyle \hbox{det}(P_{
(0,\beta_0,0,\beta_1,\dots,0,\beta_{n-1},0,\tilde \beta_n+1),-
(0,\beta_0,0,\beta_1,\dots,0,\beta_{n-1},0,\tilde \beta_n+1)}(2n+2))=1\cr
\displaystyle \hbox{det}(P_{
(0,\beta_0,0,\beta_1,\dots,0,\beta_{n-1},0,\tilde \beta_n),-
(0,\beta_0,0,\beta_1,\dots,0,\beta_{n-1},0,\tilde \beta_n)}(2n+2))=0\cr
\displaystyle \hbox{det}(P_{
(0,\beta_0,0,\beta_1,\dots,0,\beta_{n-1},0,\tilde \beta_n-1),-
(0,\beta_0,0,\beta_1,\dots,0,\beta_{n-1},0,\tilde \beta_n-1)}(2n+2))=1\ .
\end{array}$$

\ \ (ii) If $\beta=(\beta_0,\beta_1,\beta_2,\dots)$ and
$\beta'=(\beta_0',\beta_1',\beta_2',\dots)$ are two infinite sequences
of integers satisfying the assumption of assertion (i) above for all 
$n$, then there exists a unique integer $m$ such that
$\beta_i=\beta'_i$ for $i<m$ and $\beta_m=\tilde \beta_m+\epsilon$,
$\beta'_m=\tilde \beta_m-\epsilon$ with $\tilde \beta_m$ as in
assertion (i) above and $\epsilon\in\{\pm 1\}$.}

{\bf Proof.} The determinant of the symplectic matrix
$$P_{(0,\beta_0,0,\beta_1,\dots,0,\beta_{n-1},0,x),-
(0,\beta_0,0,\beta_1,\dots,0,\beta_{n-1},0,x)}(2n+2)$$
is of the form $D(x)=(ax+b)^2$ for some suitable integers $a$ and $b$
(which are well defined up to multiplication by $-1$).

It is easy to see that it is enough to show that $a=\pm 1$ in order
to prove the Theorem (the integer $\tilde \beta_m$ equals then
$-ab$ and is even by a consideration $\pmod 2$). 
This is of course equivalent to showing that the polynomial
$D(x)$ has degree $2$ and leading term $1$.

Consider now the symplectic matrix $M$ of order $2n+2$ defined as follows:
The entries of $M$ except the last row and column are given by the odd-order
(and hence degenerate) symplectic matrix
$$P_{(0,\beta_0,0,\beta_1,\dots,0,\beta_{n-1},0),-
(0,\beta_0,0,\beta_1,\dots,0,\beta_{n-1},0)}(2n+1)\ .$$
The last row (which determines by antisymetry the last column)
of $M$ is given by
$$(1,1,1,\dots,1,1,0)\ .$$
It is obvious to check that $\hbox{det}(M)$ is the coefficient of $x^2$
in the polynomial $D(x)$ introduced above.

Subtract now row number $2n-1$ from row number $2n$ of $M$ (with
rows and columns of $M$ indexed from $0$ to $2n+1$), subtract then
row number $2n-2$ from row number $2n-1$ etc. until subtracting row number 
$0$ from row number $1$. Do the same operations on columns thus producing
a symplectic matrix $\tilde M$ which is equivalent to $M$ and whose last
row is given by
$(1,0,0,\dots,0,0)$.
The determinant of $M$ equals hence the determinant of the submatrix of
$\tilde M$ obtained by deleting the first and last rows and columns 
in $\tilde M$. This submatrix is given by
$$P_{(0,\beta_0,0,\beta_1,\dots,0,\beta_{n-1}),-
(0,\beta_0,0,\beta_1,\dots,0,\beta_{n-1})}(2n)$$
thus showing that $\hbox{det}(M)=1=a^2$.\hfill QED
 
The set of sequences
$$\{\alpha=(0,\beta_0,0,\beta_1,0,\beta_2,\dots)\ \vert\ 
\hbox{det}(P_{\alpha,-\alpha}(2n)=1,\ n=1,2,3,\dots\}$$
associated
to unimodular symplectic matrices $P_{\alpha,-\alpha}(2n)$
consists hence of integral sequences and has the structure of a tree. We call 
this tree the {\it even symplectic unimodular tree}.

The beginning of this tree is shown below and is to be understood as follows:

Column $i$ displays the integer $\tilde \beta_i$ of the Theorem. 
Indices indicate if $\beta_i=\tilde \beta_i+1$ or
$\tilde \beta_i-1$. Hence the row
$$\begin{array}{rrrrrr}
0_{+1}&0_{+1}&0_{-1}&-8_{+1}&68_{+1}&434748_{\pm}\end{array}$$
corresponds for instance to the 
sequence 
$$(1,1,-1,-7,69)$$ 
implying $\tilde \beta_5=434748$ (the sequence $(1,1,-1,-7,69)$ can hence
be extended either to $(1,1,-1,-7,69,434749)$ or to
$(1,1,-1,-7,69,434747)$).

We have only displayed sequences starting with $1$ since sequences starting 
with $-1$ are obtained by a global sign change.

{\bf Table 5.1.3.} (Part of the even symplectic unimodular tree).
$$\begin{array}{rrrrrr}
0_{+1}&0_{+1}&0_{+1}&0_{+1}&0_{+1}&0_{\pm}\cr
0_{+1}&0_{+1}&0_{+1}&0_{+1}&0_{-1}&-100_{\pm}\cr
0_{+1}&0_{+1}&0_{+1}&0_{-1}&-42_{+1}&32658_{\pm}\cr
0_{+1}&0_{+1}&0_{+1}&0_{-1}&-42_{-1}&-39754_{\pm}\cr
0_{+1}&0_{+1}&0_{-1}&-8_{+1}&68_{+1}&434748_{\pm}\cr
0_{+1}&0_{+1}&0_{-1}&-8_{+1}&68_{-1}&-400344_{\pm}\cr
0_{+1}&0_{+1}&0_{-1}&-8_{-1}&-254_{+1}&12922350_{\pm}\cr
0_{+1}&0_{+1}&0_{-1}&-8_{-1}&-254_{-1}&-13258926_{\pm}\cr
0_{+1}&0_{-1}&2_{+1}&0_{+1}&240_{+1}&13257990_{\pm}\cr
0_{+1}&0_{-1}&2_{+1}&0_{+1}&240_{-1}&-12923278_{\pm}\cr
0_{+1}&0_{-1}&2_{+1}&0_{-1}&-74_{+1}&400664_{\pm}\cr
0_{+1}&0_{-1}&2_{+1}&0_{-1}&-74_{-1}&-434420_{\pm}\cr
0_{+1}&0_{-1}&2_{-1}&0_{+1}&36_{+1}&39594_{\pm}\cr
0_{+1}&0_{-1}&2_{-1}&0_{+1}&36_{-1}&-32810_{\pm}\cr
0_{+1}&0_{-1}&2_{-1}&0_{-1}&2_{+1}&92_{\pm}\cr
0_{+1}&0_{-1}&2_{-1}&0_{-1}&2_{-1}&0_{\pm}\cr
\end{array}$$

\section{A \lq\lq sympletric'' tree?}

The construction of a generalized Pascal triangle $P_{\alpha,\beta}(\infty)$ needs 
two sequences $\alpha=(\alpha_0,\alpha_1,\dots)$ and $\beta=(\beta_0,\beta_1,\dots)$.
Starting with only one sequence $\alpha=(\alpha_0,\alpha_1,\dots)$ and considering 
$P_{\alpha,\alpha}(\infty)$ we get generalized symmetric Pascal triangles
and considering $P_{\alpha,-\alpha}(\infty)$ we get generalized symplectic 
Pascal triangles. Since the sequence $\tilde \alpha=(\tilde \alpha_0=\alpha_0,
\tilde \alpha_1=-\alpha_1,\tilde \alpha_2=\alpha_2,\dots,\tilde \alpha_i=
(-1)^i\alpha_i,\dots)$ is half-way
between $\alpha$ and $-\alpha$, we call the generalized Pascal triangle 
$P_{\alpha,\tilde \alpha}(\infty)$
the {\it generalized \lq\lq sympletric'' Pascal triangle}.

The two sequences
$$\begin{array}{ll}
\alpha=(0,1,1,2,3,5,8,13,\dots)\qquad &\hbox{Fibonacci}\cr
\alpha=(0,1,1,0,-1,-1,0,1,1,0,-1,-1,\dots)\qquad &6-\hbox{periodic}\end{array}$$
and the associated sequences $\tilde \alpha$ satisfy all linear recursions of 
order 2. Theorem 3.1  and a computation of the first few values
show that both sequences $\hbox{det}(P_{\alpha,\tilde \alpha}(n))$ equal
$0,1,2,2^2,2^3,\dots,2^{n-2},\dots$.
All the following finite sequences yield matrices $P_{\alpha\tilde\alpha}(n)$
with determinants 
$0,1,2,4,8,16,32,64$ (for $n=1,2,3,\dots$) too:
$$\begin{array}{rrrrrrrr}
0&1&1&2&5&13&34&85\pm 4\cr
0&1&1&2&5&13&28&79\pm 4\cr
0&1&1&2&5&9&20&77\pm 38\cr
0&1&1&2&5&9&2&-193\pm 110\cr
0&1&1&2&3&5&10&19\pm 6\cr
0&1&1&2&3&5&8&12\pm 1\cr
0&1&1&2&3&3&10&-3\pm 6\cr
0&1&1&2&3&3&-2&9\pm 30\cr
0&1&1&0&1&3&4&-31\pm 38\cr
0&1&1&0&1&3&-14&167\pm 110\cr
0&1&1&0&1&-1&2&1\pm 4\cr
0&1&1&0&1&-1&-4&-17\pm 4\cr
0&1&1&0&-1&1&10&33\pm 6\cr
0&1&1&0&-1&1&-2&-27\pm 30\cr
0&1&1&0&-1&-1&2&3\pm 6\cr
0&1&1&0&-1&-1&0&2\pm 1\cr
\end{array}$$
{\bf Problem 6.1.} Has the set of all infinite integral
sequences $\alpha=(0,1,1,\alpha_3,\dots)$
such that $\hbox{det}(P_{\alpha,\tilde \alpha}(n))=
(0,1,2,4,\dots,2^{n-2},\dots)_{n=1,2,\dots}$ the structure of a tree
(ie. can every finite such sequence of length at least $3$
be extended by one next term in exactly two ways)? 

\section{Periodic matrices}

In this section we are interested in matrices coming from a kind of
\lq\lq periodic convolution with compact support on ${\bf N}$''.  

We say that an infinite matrix $A$ with coefficients $a_{i,j},\ 0\leq i,j$
is {\it $(s,t)-$bounded} ($s,t\in{\bf N}$) if $a_{i,j}=0$ for $(j-i)\not\in [-s,t]$.

We call a matrix with coefficients $a_{i,j},\ 0\leq i,j$ 
{\it $p-$periodic} if $a_{i,j}=a_{i-p,j-p}$ for $i,j\geq p$.

An infinite matrix $P$ with coefficients $p_{i,j},\ 0\leq i,j$ is a 
{\it finite perturbation} if it has only a finite number of non-zero
coefficients.

As before, given an infinite matrix $M$ with coefficients $m_{i,j},\ 0\leq i,j$
we denote by $M(n)$ the matrix with coefficients $m_{i,j},\ 0\leq i,j<n$
obtained by erasing all but the first
$n$ rows and columns of $M$. 

{\bf Theorem 7.1.} {\sl Let $A=\tilde A+P$ be a matrix where $\tilde A$ is
 a $p-$periodic $(s,t)-$bounded  matrix and where
$P$ is a finite perturbation.
Then there exist constants $N,d\leq {s+t\choose s},C_1,\dots,C_d$ such that
$$\hbox{det}(A(n))=\sum_{i=1}^d C_i\ \hbox{det}(A(n-ip))$$
for $n>N$.}

We will prove the theorem for $p=1$, $s=t=2$ and then describe the necessary
modifications in the general case.

{\bf Proof in the case $p=1,\ s=t=2$.} Suppose $n$ huge. The matrix $A(n)$
has then the form
$$A(n)=\left(\begin{array}{ccccc}
\ddots&&&&\vdots\cr
&c&d&e&0\cr
&b&c&d&e\cr
&a&b&c&d\cr
\dots&0&a&b&c\end{array}\right)\ .$$
Developping the determinant possibly several times along the last row one
gets only matrices of the following six types
$$T_1=\left(\begin{array}{cccc}
\ddots&d&e&0\cr
&c&d&e\cr
&b&c&d\cr
&a&b&c\end{array}\right)
\quad 
T_2=\left(\begin{array}{cccc}
\ddots&d&e&0\cr
&c&d&0\cr
&b&c&e\cr
&a&b&d\end{array}\right)
\quad 
T_3=\left(\begin{array}{cccc}
\ddots&d&0&0\cr
&c&e&0\cr
&b&d&e\cr
&a&c&d\end{array}\right)$$
$$
T_4=\left(\begin{array}{cccc}
\ddots&d&e&0\cr
&c&d&0\cr
&b&c&0\cr
&a&b&e\end{array}\right)
\quad 
T_5=\left(\begin{array}{cccc}
\ddots&d&0&0\cr
&c&e&0\cr
&b&d&0\cr
&a&c&e\end{array}\right)
\quad 
T_6=\left(\begin{array}{cccc}
\ddots&d&0&0\cr
&c&0&0\cr
&b&e&0\cr
&a&d&e\end{array}\right)
$$
and writing $t_i(m)=\hbox{det}(T_i(m))$ we have the identity
$$\left(\begin{array}{c}
t_1(m)\cr t_2(m)\cr t_3(m)\cr t_4(m)\cr t_5(m)\cr t_6(m)\end{array}\right)
=\left(\begin{array}{cccccc}
c&-b&a&0&0&0\cr
d&0&0&-b&a&0\cr
0&d&0&-c&0&a\cr
e&0&0&0&0&0\cr
0&e&0&0&0&0\cr
0&0&0&e&0&0\end{array}\right)
\left(\begin{array}{c}
t_1(m-1)\cr t_2(m-1)\cr t_3(m-1)\cr t_4(m-1)\cr t_5(m-1)\cr t_6(m-1)
\end{array}\right)$$
for $m$ huge enough. Writing $R$ the above $6\times 6$ matrix
relating $t_i(m)$ to $t_j(m-1)$ we have
$t(n)=R^{n-N}t(N)$ for $n\geq N$ huge enough and for $t(m)$ the vector
with coordinates $t_1(m),\dots,t_6(m)$. Choosing a basis of a Jordan normal form
of $R$ and expressing the vector $t(N)$ with respect to this basis shows now that 
the determinants $t_i(n)$ (and hence $\hbox{det}(A(n))=t_1(n)$)
satisfy for $n>N$
a linear recursion with characteristic polynomial dividing
$$\hbox{det}(z\ \hbox{Id}_6-R)\ .$$

{\bf Proof of the general case.} Let us first suppose $p=1$. There are then
${s+t\choose s}$ (count the possibilities for the highest non-zero entry
in the last $s$ columns) different possible types $T_i$
obtained by developping the determinant $\hbox{det}(A(n))$ for huge $n$
several times along the last row and one gets hence a square matrix $R$ of order
${s+t\choose s}$ expressing the determinants $\hbox{det}(T_i(n))$ linearly
in $\hbox{det}(T_j(n-1))$ for $n$ huge enough. This shows that the 
determinants $\hbox{det}(T_i(n))$ satisfy for $n$ huge enough
a linear recursion with
characteristic polynomial dividing the characteristic polynomial of 
the square matrix $R$.

If $p>1$, develop the determinant of $\hbox{det}(A(n))$ a multiple of $p$ times
along the last row and proceed as above. One gets in this way matrices
$R_0,\dots,R_{p-1}$ according to $n\pmod p$ with identical characteristic polynomials
yielding recursion relations between $\hbox{det}(A(n))$ and $\hbox{det}(A(n-ip))$.
\hfill QED
 
\section{The diagonal construction}

Let $\gamma=(\gamma_0,\gamma_1,\gamma_2,\dots,)$ be a sequence and let 
$u_1,u_2,l_1,l_2$ be four constants. The {\it diagonal-construction}
is the (infinite) matrix $D_\gamma^{(u_1,u_2,l_1,l_2)}$ with
entries
$$\begin{array}{ll}
d_{i,i}=\gamma_i\ &0\leq i\cr
d_{i,j}=u_1d_{i,j-1}+u_2d_{i+1,j}\quad &0\leq i<j\cr
d_{i,j}=l_1d_{i-1,j}+l_2d_{i,j+1}\ &0\leq j<i\end{array}$$
and we denote by $D(n)=D_\gamma^{(u_1,u_2,l_1,l_2)}(n)$ the
$n\times n$ principal submatrix with coefficients $d_{i,j},\ 0\leq i,j<n$
obtained by considering the first $n$
rows and columns of $D_\gamma^{(u_1,u_2,l_1,l_2)}$.

The cases where $u_1u_2l_1l_2=0$ are degenerate. For instance, in the case 
$u_2=0$ one sees easily that the matrix 
$D_\gamma^{(u_1,u_2,l_1,l_2)}(n)$ has determinant
$$\gamma_0\ \prod_{j=1}^{n-1}(\gamma_j-u_1(l_1\gamma_{j-1}+l_2\gamma_j))\ .$$
The other cases are similar.

The following result shows that we loose almost nothing by assuming
$u_1=u_2=1$.

{\bf Proposition 8.1.} {\sl For $\lambda,\mu$ two invertible constants we
have
$$D_{\tilde \gamma}^{(\lambda u_1,\mu u_2,\mu^{-1}l_1,\lambda^{-1}l_2)}(n)=
\left(\frac{\lambda}{\mu}\right)^{n\choose 2}\ 
D_\gamma^{(u_1,u_2,l_1,l_2)}(n)$$
where 
$$\tilde \gamma=(\gamma_0,\frac{\lambda}{\mu}\gamma_1,
\frac{\lambda^2}{\mu^2}\gamma_2,\dots,\tilde\gamma_k=
\frac{\lambda^k}{\mu^k}\gamma_k,\dots)\ .$$}

{\bf Proof.} Check that the coefficients $\tilde d_{i,j}$ of
$D_{\tilde \gamma}^{(\lambda u_1,\mu u_2,\mu^{-1}l_1,\lambda^{-1}l_2)}(n)$
are given by $\tilde d_{i,j}=\mu^{-i}\lambda^j\ d_{i,j}$ where $d_{i,j}$
are the coefficients of $D_\gamma^{(u_1,u_2,l_1,l_2)}(n)$.
This implies the result easily.\hfill QED

{\bf Proposition 8.2.} {\sl For $n\geq 1$ the sequence
$$d(n)=\hbox{det}(D_\gamma^{(u_1,u_2,l_1,l_2)}(n))$$
associated to the geometric sequence $\gamma=(1,x,x^2,x^3,\dots)$ is
given by
$$d(n)=\left(-u_1l_1+(1-u_1l_2-u_2l_1)x-u_2l_2x^2\right)^{n-1}\ 
x^{n-1\choose 2}\ .$$}

A nice special case is given by $u_1=u_2=l_1=l_2=1$. The associated matrix
$D_\gamma(4)=D_\gamma^{(1,1,1,1)}(4)$ for example is then given by
$$\left(\begin{array}{cccc}
1&1+x&1+2x+x^2&1+3x+3x^2+x^3 \cr
1+x&x&x+x^2&x+2x^2+x^3\cr
1+2x+x^2&x+x^2&x^2&x^2+x^3\cr
1+3x+3x^2+x^3&x+2x^2+x^3 &x^2+x^3&x^3
\end{array}\right)$$
and the reader can readily check that the coefficient $d_{i,j}$ of $D_\gamma(n)$
is given by
$$d_{i,j}=x^{\hbox{min}(i,j)}(1+x)^{\vert i-j\vert}\ .$$
Proposition 8.2 shows that the determinant $\hbox{det}(D_\gamma(n))$
is given by
$$\hbox{det}(D_{(1,x,x^2,x^3,\dots)}(n))=\left(-1-x-x^2\right)^{n-1}x^{n-1\choose 2}$$
for $n\geq 1$.

Setting $x=1$ in this special case $u_1=u_2=l_1=l_2=1$, we get a 
matrix $M$ with
entries $m_{i,j}=2^{\vert i-j\vert}$ for $0\leq i,j<n$. Its determinant is 
$(-3)^{n-1}$. It is easy to show that the matrix $M_a$ of order $n$ with 
entries $m_{i,j}=a^{\vert i-j\vert}$ for $0\leq i,j<n$ 
has determinant $(1-a^2)^{n-1}$. 

A similar example is the special case $-u_1=u_2=-l_1=l_2=1$ which yields for
instance the matrix $D_\gamma(4)=D_\gamma^{(-1,1,-1,1)}(4)$ given by
$$\left(\begin{array}{cccc}
1&-1+x&1-2x+x^2&-1+3x-3x^2+x^3 \cr
-1+x&x&-x+x^2&x-2x^2+x^3\cr
1-2x+x^2&-x+x^2&x^2&-x^2+x^3\cr
-1+3x-3x^2+x^3&x-2x^2+x^3 &-x^2+x^3&x^3
\end{array}\right)$$
and the reader can readily check that the coefficient $f_{i,j}$ of $D_\gamma(n)$
is given by
$$d_{i,j}=x^{\hbox{min}(i,j)}(x-1)^{\vert i-j\vert}\ .$$

The determinant $\hbox{det}(D_{(1,x,x^2,x^3,\dots)}(n))$
is given by
$$\hbox{det}(D_{(1,x,x^2,x^3,\dots)}(n))=\left(-x^2+3x-1\right)^{n-1}
x^{n-1\choose 2}$$
for $n\geq 1$.

{\bf Proof of Proposition 8.2.} By continuity and Proposition 8.1 it is enough 
to prove the formula in the case $u_1=u_2=1$.

This implies $d_{i,j}=x^i(1+x)^{(j-i)}$ for $i\leq j$.

Subtracting $(1+x)$ times column number $(n-2)$ from column number $(n-1)$
(which is the last one) etc. until subtracting $(1+x)$ times column number $0$ 
from column number $1$ transforms the matrix $D(n)$ into a lower
triangular matrix with diagonal entries 
$$1,x-(1+x)(l_1+l_2x),x\left(x-(1+x)(l_1+l_2x) \right),\dots,
x^{n-2}(-l_1+(1-l_1-l_2)x-l_2x^2)\ .$$ 

{\bf Theorem 8.3.} {\sl Let $\gamma=(\gamma_0,\dots,\gamma_{p-1},\gamma_0,
\dots,\gamma_{p-1},\dots)$ be a $p-$periodic sequence and let
$$d(n)=\hbox{det}(D_\gamma^{(u_1,u_2,l_1,l_2)}(n))$$
be the determinants of the associated matrices (for fixed 
$(u_1,u_2,l_1,l_2)$).

Then there exist an integer $d$ and constants $C_1,\dots
C_{d}$ such that
$$d(n)=\sum_{i=1}^{d}C_i\ d(n-ip)$$
for all $n$ huge enough.}

{\bf Remark 8.4.} Generically, the coefficients $C_i$ seem to display the
symmetry
$$C_{d-i}=\rho^{(d-2i)/2}\ C_i$$
(with $C_0=-1$) 
for some constant $\rho$ which seems to be polynomial in
$\gamma_0,\dots,\gamma_{p-1},u_1,u_2,l_1,l_2$.

{\bf Proof of Theorem 8.3.} For $k\geq p$ add to the $k-$th row a linear 
combination (with coefficients depending only on $(u_1,u_2,l_1,l_2)$) of 
rows $k-1,k-2,\dots,k-p$ such that $d_{k,k-i}=0$ for $i\geq p$. Do the
analogous operation on columns and apply Theorem 7.1 to the resulting 
matrices. \hfill QED

I thank D. Fux, P. de la Harpe and especially C. Krattenthaler for many usefull 
remarks and important improvements over a preliminary version.

\vskip1cm

{\bf Bibliography}

{\bf [GHJ]} Goodman,P. de la Harpe, V.F.R. Jones, {\it Coxeter graphs and towers of 
algebras}, Springer (1989).

{\bf [GV]} I. Gessel, G. Viennot, {\it Binomial Determinants, Paths,
and Hook Length Formulae}, Adv. Math. {\bf 58} (1985), 300-321.

{\bf [IS]} Integer-sequences, http://www.research.att.com/~njas/sequences/index.html

{\bf [K1]} C. Krattenthaler, {\it Advanced Determinant Calculus}, 
S\'eminaire Lotharingien Combin. 42 ("The Andrews Festschrift") (1999),
  Article B42q, 67 pp. 

{\bf [K2]} C. Krattenthaler, {\it personal communication}, september 2001.

{\bf [MW]} M.L. Mehta, R. Wang, {\it Calculation of a certain determinant}, 
Commun. Math. Phys. 214, No.1, 227-232 (2000). 

\vskip1cm

Roland Bacher

INSTITUT FOURIER

Laboratoire de Math\'ematiques

UMR 5582 (UJF-CNRS)

BP 74

38402 St MARTIN D'H\`ERES Cedex (France)
 
e-mail: Roland.Bacher@ujf-grenoble.fr

\end{document}